\documentclass[a4paper,12pt]{article}
\usepackage{amsmath}
\usepackage{amssymb}
\usepackage{amsthm} 

\usepackage{tikz}
\usetikzlibrary{arrows.meta,positioning}

\theoremstyle{plain}
\newtheorem{theorem}{Theorem}

\newtheorem{claim}[theorem]{Claim}
\newtheorem{corollary}[theorem]{Corollary}

\theoremstyle{definition}
\newtheorem{definition}[theorem]{Definition}
\newtheorem{remark}[theorem]{Remark}

\author{Hiroki Kodama}
\title
{
Characterization of cactus-expandable digraphs via doubly bidirectionally connected pairs
}%
\date{Jun.\ 3rd, 2026}

\begin{document}

\maketitle

\begin{abstract}
Azuma et al.\ showed that a strongly connected digraph without a doubly bidirectionally connected pair is cactus-expandable.
We prove the converse: if a digraph has a doubly bidirectionally connected pair, then no expansion of it is a cactus digraph.
Combined with the theorem of Azuma et al., this yields a characterization of strongly connected cactus-expandable digraphs.
\end{abstract}

The main result of this paper is the following theorem.

\begin{theorem}\label{main}
If a digraph $G$ has a doubly bidirectionally connected pair $(p,q)$ and $\varphi \colon G' \to G$ is an expansion, then
$G'$ is not a cactus digraph.
\end{theorem}

Combining Theorem \ref{main} with the theorem of Azuma et al. \cite{Az} immediately yields the following characterization.

\begin{corollary}
Suppose a digraph $G$ is strongly connected.
Then the following are equivalent:
\begin{itemize}
\item There exists an expansion $\varphi \colon G' \to G$ such that $G'$ is a cactus digraph.
\item $G$ does not have a dbcp.
\end{itemize}
\end{corollary}

At first glance, Theorem \ref{main} might appear to follow directly from the structural rigidity of expansions. 
However, the obstruction caused by a doubly bidirectionally connected pair is more subtle, as it persists under all possible expansions.

We remark that our notion of a cactus digraph differs from the standard notion of a cactus graph in undirected graph theory.

\bigskip

A digraph, or directed graph, is an ordered pair $G=(V,A)$,
where $V$ is a finite set of vertices and $A$ is a finite set of ordered pairs (called arcs) of vertices of $V$, in other words, $A\subset V \times V$.
In this paper we do not consider loops, i.e.\ $(v,v)\not\in A$ for any $v\in V$.

For an arc $a=(x,y)$, we use the notation that $s(a)=x$ and $t(a)=y$. We also say that
$a$ is an arc from $x$ to $y$.
For a vertex $v \in V$, 
the indegree $\delta^-(v)$ is the number of arcs to $v$ and
the outdegree $\delta^+(v)$ is the number of arcs from $v$.

A path of length $n$ from $x$ to $y$ is a sequence of $n+1$ vertices $(x_0,x_1,\dots,x_{n-1},x_n)$ such that
$x_0=x$, $x_n=y$ and $(x_{i-1},x_i)\in A$ for $1 \leq i \leq n$.
The vertex $x_{n-1}$ is called the preterminal vertex of the path. 
A path is called simple if $x_0, \dots, x_{n-1}, x_n$ are distinct.

A cycle is a path with $x_0 = x_n$.
A cycle is called simple if $x_0, \dots, x_{n-1}$ are distinct. Note that a simple cycle cannot be a simple path.

A digraph $G=(V,A)$ is said to be strongly connected if for any vertices $x,y \in V$ there exists a path from $x$ to $y$.

\begin{definition}[morphism]
Suppose $G=(V,A)$ and $G'=(V',A')$ are digraphs and 
$\varphi=(\varphi_V, \varphi_A)$ is a pair of maps s.t.\ 
$\varphi_V \colon V' \to V$ and
$\varphi_A \colon A' \to A$.
$\varphi$ is said to be a digraph morphism from $G'$ to $G$ if
$\varphi_V$ and $\varphi_A$ are compatible, 
i.e.\ for any $(v', w') \in A'$,
$\varphi_A((v', w'))=(\varphi_V(v'), \varphi_V(w'))$.
To denote the digraph morphism, we write $\varphi \colon G' \to G$.
\end{definition}

The idea of expansions of digraphs is introduced by Azuma et al.\ in their paper \cite{Az}, where they analyzed 
Boolean dynamical systems characterized by digraphs.

\begin{definition}[expansion] 
Suppose $\varphi \colon G' \to G$ is a digraph morphism from $G'=(V',A')$ to $G=(V,A)$.
An ordered pair $(G',\varphi)$ is said to be an expansion of a digraph $G$ if the following condition holds:
\begin{itemize}
\item $\varphi_V \colon V' \to V$ is a surjective map.
\item $\varphi_A \colon A' \to A$ is a surjective map.
\item If $(x,y) \in A$, $y' \in V'$ and $y = \varphi_V(y')$, there exists a unique arc $(x',y') \in A'$ such that $\varphi_V(x') = x$.
\end{itemize}
\end{definition}

For $v' \in V'$ and $a' \in A'$, we often write $\varphi(v')$ and $\varphi(a')$ instead of
$\varphi_V(v')$ and $\varphi_A(a')$. We also say that $\varphi \colon G' \to G$ is an expansion.
Note that the last condition implies that the indegree $\delta^-$ is preserved under expansion, i.e., if $\varphi_V(v') = v$, then $\delta^-(v') = \delta^-(v)$, 
but this does not hold for the outdegree $\delta^+$.

Here is an example.

\begin{figure}[htbp]
\centering
\begin{tikzpicture}[
    scale=0.8,
    transform shape,
    >=Stealth,
    every node/.style={circle,draw,inner sep=1.5pt,minimum size=16pt,font=\small},
    edge/.style={->,thick}
]

\node (v1) at (0,0)   {$v_1$};
\node (v2) at (-1.8,-1.3) {$v_2$};
\node (v3) at (1.8,-1.3) {$v_3$};

\node (u1) at (7.2,0) {$u_1$};
\node (u2) at (5.4,-1.3) {$u_2$};
\node (u3) at (9.0,-1.3) {$u_3$};
\node (u4) at (10.8,0) {$u_4$};

\draw[edge] (v1) -- (v2);
\draw[edge] (v2) -- (v3);
\draw[edge] (v3) -- (v1);
\draw[edge] (v1) ..controls+(1.4,-0.25).. (v3);

\draw[edge] (u1) -- (u2);
\draw[edge] (u2) -- (u3);
\draw[edge] (u3) -- (u1);
\draw[edge] (u3) -- (u4);
\draw[edge] (u4) ..controls+(-0.4,-1.05).. (u3);

\end{tikzpicture}
\caption{An example of an expansion.}\label{fig:expansion}
\end{figure}

Let $\varphi$ be a digraph morphism from the graph on the right-hand side of Figure \ref{fig:expansion} to the graph on the left-hand side, 
defined by
\[
\varphi_V(u_1)=\varphi_V(u_4)=v_1,\quad \varphi_V(u_2)=v_2,\quad \varphi_V(u_3)=v_3,
\]
and let $\varphi_A$ be defined so as to be compatible with $\varphi_V$. Then $\varphi$ defines an expansion.

\begin{definition}[cactus digraph]
A cactus digraph is a strongly connected digraph in which each arc is contained in exactly one simple cycle.

We remark again that our notion of a cactus digraph differs from the standard notion of a cactus graph in undirected graph theory.
\end{definition}

In particular, in a cactus digraph, each arc belongs to one simple directed cycle and
any two simple directed cycles share at most one vertex.


\begin{figure}[htbp]
\centering
\begin{tikzpicture}[
    scale=0.8,
    transform shape,
    >=Stealth,
    every node/.style={circle,draw,inner sep=1.5pt,minimum size=16pt,font=\small},
    edge/.style={->,thick}
]

\node (v1) at (0,0)   {$v_1$};
\node (v2) at (1.8,1.3) {$v_2$};
\node (v3) at (3.6,0) {$v_3$};
\node (v4) at (1.8,-1.3) {$v_4$};

\node (v5) at (5.4,1.3) {$v_5$};
\node (v6) at (7.2,0) {$v_6$};

\node (v7) at (3.6,-2.6) {$v_7$};
\node (v8) at (7.2,-2.6) {$v_8$};

\node (v9)  at (-1.8,1.3) {$v_9$};
\node (v10) at (-3.6,0) {$v_{10}$};

\draw[edge] (v1) -- (v2);
\draw[edge] (v2) -- (v3);
\draw[edge] (v3) -- (v4);
\draw[edge] (v4) -- (v1);

\draw[edge] (v3) -- (v5);
\draw[edge] (v5) -- (v6);
\draw[edge] (v6) -- (v3);

\draw[edge] (v3) -- (v7);
\draw[edge] (v7) -- (v8);
\draw[edge] (v8) -- (v3);

\draw[edge] (v1) -- (v9);
\draw[edge] (v9) -- (v10);
\draw[edge] (v10) -- (v1);

\end{tikzpicture}
\caption{An example of a cactus digraph.}\label{fig:cactus}
\end{figure}


Figure \ref{fig:cactus} and right-hand side of Figure \ref{fig:expansion} are examples of cactus digraphs.

\begin{remark}
We call the digraph with single vertex and no arc to be a cactus digraph.
\end{remark}

\begin{claim}
A strongly connected digraph is a cactus digraph if and only if two distinct simple cycles share at most one vertex.
\end{claim}
\begin{proof}
If an arc $(x,y)$ is contained by two simple cycles, two vertices $x$ and $y$ are also shared by those two cycles.
On the other hand, suppose two distinct simple cycles $S_1$ and $S_2$ share two vertices $x$ and $y$,
then $S_1=p_1 \cup q_1$ and $S_2=p_2 \cup q_2$ where $p_i$ is a simple path from $x$ to $y$ and $q_i$ from $y$ to $x$
($i=1,2$). 
If the preterminal points of $p_1$ and $p_2$ coincide, let that vertex be $y'$. Then the arc $(y',y)$ is shared by
 $S_1$ and $S_2$, that violates the definition of a cactus digraph.
So we can assume that they are different points and $p_1 \neq p_2$.
Then $S_1=p_1 \cup q_1$ and $S_2'=p_2 \cup q_1$ are distinct cycles.
If $S_2'$ is simple, then any arc in $q_1$ is shared by $S_1$ and $S_2'$. 

If $S_2'$ is not simple, let $q_1=(y=y_0,y_1,\dots,x)$ and $i$ be the smallest $i>0$ s.t.\ $y_i\in p_1$.
Let $p_1'$ be the sub path of $p_1$ from $z=y_i$ to $y$ and
$q_1'$  the sub path of $q_1$ from $y$ to $z$.
Then $S_1$ and $S_2''=p_1'\cup q_1'$ are distinct simple cycles, and
any arc in $q_1'$ is shared by $S_1$ and $S_2''$.
\end{proof}

\begin{definition}[connecting point]
A vertex in a cactus digraph shared by two or more distinct simple cycles is called a connecting point.
\end{definition}

A cactus digraph is balanced, i.e.\ for any vertex $v \in V$, $\delta^-(v)=\delta^+(v)$.
A connecting point is a vertex $v$ s.t.\ $\delta^-(v)=\delta^+(v) \geq 2$.

If $G_1$ and $G_2$ are sub cactus digraphs of a cactus digraph $G$, then the intersection $G_1 \cap G_2$ is
either sub cactus digraph of $G$, a single vertex, or an empty set. For $\tilde{V} \subset V$, consider all 
sub cactus digraphs of $G$ that contain $\tilde{V}$ and take those intersection $\tilde{G}$. Then
$\tilde{G}$ is the minimum cactus subgraph that contains $\tilde{V}$. Remark that in case $\tilde{V}$
is a single vertex set $\{v\}$, $\tilde{G}=\{v\}$.

\bigskip

We would like to define a preorder for a rooted cactus digraph, but before that
we review the definition of preorder.

\begin{definition}[preorder]
A binary relation $\precsim$ of a set $P$ is called a preorder if it is
\begin{itemize}
\item reflexive, i.e.\ $x \precsim x$ for any $x \in P$, and
\item transitive, i.e.\ if $x \precsim y$ and $y \precsim z$ then $x \precsim z$ for any $x,y,z \in P$.
\end{itemize}
\end{definition}

In other words, preorder is a partial order minus antisymmetric law.
If $x \precsim y$ and $y \precsim x$, we write $x \sim y$, but $x$ may differ from $y$.
In case $x \precsim y$ and $y \not\precsim x$, we write $x \prec y$.

Fix a root point $r \in V$, then we can define the preorder for the rooted cactus digraph $(G,r)$
in the following way.

\begin{definition}[preorder for the rooted cactus digraph]
For $v \in V\setminus \{r\}$, we define $C(v)$ as the minimum sub cactus digraph of $G$ that contains
$\{r,v\}$. For $r$, $C(r) = \{r\}$. 
Then preorder in $V$ is defined by $v \precsim w$ if and only if $C(v) \subseteq C(w)$.
\end{definition}

Here are some properties of the preorder for the rooted cactus digraph.

The root $r$ is the strongly minimum element, i.e.\ for any $v\in V\setminus \{r\}$ $r \prec v$.
On the other hand, since $V$ is a finite set, there exists a weakly maximal element $m \in V$, i.e.\ for any $v\in V$ $m \not\prec v$.

\begin{claim}
If $S=(V_S, A_S)$ is a simple cycle of the cactus digraph $G$, then there exists the strongly minimum element 
$v_S$ among $V_S$. $v_S=r$ if $S$ contains the root $r$ and $v_S$ is the nearest connecting point to $r$ otherwise.
We call $v_S$ the minimum point of $S$.
For $v,w \in V_S \setminus \{v_S\}$, $v\sim w$.
\end{claim}

\begin{proof}
If $r \in V_S$, the claim is trivial.
If $r \not\in V_S$, consider 
the minimum sub cactus digraph $C$ that contains $V_S \cup \{r\}$. Then $S$ is a cycle in $C$ and its unique connecting point 
in $C$ is $v_S$.
\end{proof}

\begin{claim}\label{connectingpoint}
Suppose $c \neq r$ is a connecting point of the cactus digraph $G$ s.t.\ $\delta^-(c)=\delta^+(c) = n \geq 2$.
Incoming arcs are $(x_1,c),\dots,(x_n,c)$ and 
outgoing arcs are $(c,y_1),\dots,(c,y_n)$.
Then $C(c)$ contains exactly one of $x_i$'s and one of $y_j$'s.
We can assume they are $x_1$ and $y_1$, then it follows that
$x_1,y_1 \precsim c \prec x_i,y_j$ for $2\leq i,j \leq n$.
\end{claim}

\begin{proof}
The point $c$ is not a connecting point of subcactus digraph $C(c)$, 
because if it is a connecting point, $C(c)$ can be decomposed into $C(c)= C_1 \cup \cdots \cup C_l$, 
where $l \ge 2$, each $C_s$ is a cactus subdigraph containing $c$, and $C_s \cap C_t = \{c\}$ for $s\neq t$.
There is a $C_s$ that contains the root $r$. Since $C_s$ contains both $r$ and $c$, and is strictly smaller than $C(c)$, 
this contradicts the minimality of $C(c)$.
Therefore, 
$c$ has one incoming arc $(x_i,c)$ and
one outcoming arc $(c,y_j)$ in $C(c)$.
\end{proof}


\begin{claim}\label{dipped}
Suppose $(v_0, v_1, \dots, v_n)$ is a simple path of length $n$ in the
cactus digraph $G$. Then the sequence $v_0, v_1, \dots, v_n$ is single dipped. 
i.e.\ there exists $k \in [0,n]$ s.t.\ $v_0 \succsim \cdots \succsim  v_k \precsim \cdots \precsim v_n$.
\end{claim}

\begin{proof}
Consider the subcactus digraph $C(v_0)$. If $v_0, v_1, \dots, v_n \in C(v_0)$, then
$v_0 \succsim \cdots \succsim  v_n$ so the sequence is single dipped for $k=n$.
If not, suppose $k'$ is the smallest number s.t.\ $v_{k'} \not \in C(v_0)$ and let $k= k'-1$.
Then $v_0 \succsim \cdots \succsim  v_k \prec v_{k+1} \precsim \cdots \precsim v_n$.
\end{proof}

\bigskip

The idea of doubly bidirectionally connected pair is introduced by Azuma et al.\ \cite{Az}
as a criterion for cactus expandability.

\begin{definition}[doubly bidirectionally connected pairs]
We say a pair of vertices $(p,q)\in V \times V$ of digraph $G=(V,A)$ is a doubly bidirectionally connected pair (dbcp for short) if:
\begin{itemize}
\item $p \neq q$, 
\item there exist two simple paths from $p$ to $q$ with different preterminal vertices, and
\item there exist two simple paths from $q$ to $p$ with different preterminal vertices.
\end{itemize}
\end{definition}

Then the following theorem \cite{Az} is  shown.

\begin{theorem}\label{Az}
If a digraph $G$ is strongly connected and does not have a dbcp, then
there exists an expansion $\varphi \colon G' \to G$ such that $G'$ is a cactus digraph.
\end{theorem}

Here we are going to prove the converse.


Note that strong connectivity is inherited via surjective morphisms. Therefore,
we may assume that $G$ is strongly connected, since otherwise no expansion of $G$ can be a cactus digraph.

%
%

\bigskip

Now we are going to prove the theorem \ref{main} by contradiction.

The idea of the proof is to prove that the presence of a doubly bidirectionally connected pair forces the existence of two distinct directed cycles sharing more than one vertex in any expansion, which contradicts the cactus condition.

Assume that $G$ has a dbcp $(p,q)$, $G'$ is a cactus digraph and
$\varphi \colon G' \to G$ is an expansion.

Suppose $X=\varphi_V^{-1}(\{p,q\})$ and take a root point $r$ from $X$.
Then $\precsim$ also defines the preorder in $X$.
The point $r$ is the minimum point among $X$. Take a weakly maximal point $x$ among $X$.
Since $X$ has at least two points, $x \neq r$.

We can assume that $\varphi_V(x)=q$ without loss of generality. 
$\delta^-(x) = \delta^-(q) \geq 2$, so $x$ is a connecting point of the cactus digraph $G'$.

Since there are two paths from $p$ to $q$ with different preterminal vertices, 
there are two paths from $y$ to $x$ and from $z$ to $x$ with different preterminal vertices 
s.t.\ $\varphi_V(y) = \varphi_V(z) = p$. 
Suppose vertices of these paths are 
$y=y_0,\dots,y_{m-1},y_m=x$ and
$z=z_0,\dots,z_{n-1},z_n=x$. The preterminal vertices are different, so $y_{m-1} \neq z_{n-1}$

Then from claim \ref{connectingpoint}
at least one of $y_{m-1}, z_{n-1}$ is strictly bigger than $x$. 
Suppose $z_{n-1}\succ x$ without loss of generality.
Since paths are single dipped (claim \ref{dipped}), 
$$
z=z_0 \succsim z_1 \succsim \cdots \succsim z_{n-1}\succ x.
$$
Therefore, $z \succ x$. But $z\in X$ and it contradicts that $x$ is a weakly maximal element among $X$. q.e.a.

\bigskip

\section*{Acknowledgement}
I would like to express my gratitude for the meaningful discussions with Shun-ichi Aszuma. Additionally, 
I appreciate Ryosuke Iritani, Christy Kelly, and Yukio Koriyama for their contributions and support.

\section*{AI usage disclosure}
The author used ChatGPT-5.5 to find typos and gramatical errors.
It was not used to generate mathematical results or proofs.
All AI-assisted outputs were reviewed, edited, and verified by the author, who takes full responsibility for the content
of the manuscript.

\noindent
Address: \\
International Institute for Sustainability with Knotted Chiral Meta Matter, Hiroshima University, 1-3-1 Kagamiyama, Higashi-Hiroshima, Hiroshima 739-8531 JAPAN \\
Center for Interdisciplinary Theoretical and Mathematical Sciences, RIKEN, 2-1 Hirosawa, Wako-shi, Saitama 351-0198 JAPAN

\end{document}